# On a series of Furdui and Qin and some related integrals

(Solution to Problem H-691 from the Fibonacci Quarterly with Additions)

by **Khristo N. Boyadzhiev**
Department of Mathematics and Statistics, Ohio Northern University, Ada, OH 45810
k-boyadzhiev@onu.edu

## Abstract

In this note we present the solution of Problem H-691 (published in [3]) with more details and some corrections. The solution involves three nontrivial integrals whose values were given in [3] without proof (two of them, however, with references). These integrals are

$$I_1 = \int_0^1 \frac{\ln(1+x^2)}{1+x^2} dx, \quad I_2 = \int_0^1 \frac{\ln(1+x^2)}{1+x} dx, \quad I_3 = \int_0^1 \frac{\arctan x}{1+x} dx,$$

Here we provide evaluations of these integrals in three appendices.

## Introduction and solution

Ovidiu Furdui and Huizeng Qin proposed the following problem (The Fibonacci Quarterly, Vol. 47, No. 33, August 2009/2010).

**Problem H-691**. Evaluate

$$\sigma = \sum_{n=1}^{\infty} (-1)^n \left( \ln 2 - \frac{1}{n+1} - \frac{1}{n+2} - \ldots - \frac{1}{2n} \right)^2$$

The solution was published in [3]. The value of this sum involves the Catalan constant $G$, and also $\pi$ and $\log 2$:

$$\sigma = \frac{G}{2} + \frac{\pi^2}{48} - \frac{7}{8}(\ln 2)^2 - \frac{\pi}{8} \ln 2 \ . \tag{1}$$

**Solution.** First we use a well-known identity (see, for example, [6]),

$$\frac{1}{n+1} + \frac{1}{n+2} + \ldots + \frac{1}{2n} = \sum_{k=1}^{2n} \frac{(-1)^{k-1}}{k} \ . \qquad (2)$$

At the same time,

$$\ln 2 = \sum_{k=1}^{\infty} \frac{(-1)^{k-1}}{k} \ . \qquad (3)$$

Thus

$$\ln 2 - \frac{1}{n+1} - \frac{1}{n+2} - \ldots - \frac{1}{2n} = \sum_{k=2n+1}^{\infty} \frac{(-1)^{k-1}}{k} = \int_0^1 \frac{x^{2n} dx}{1+x} . \qquad (4)$$

The last equality is easy to establish by expanding $\frac{1}{1+x}$ in power series and integrating termwise. Next we write

$$\sigma = \sum_{n=1}^{\infty} (-1)^n \left( \int_0^1 \frac{x^{2n} dx}{1+x} \right)^2 = \sum_{n=1}^{\infty} (-1)^n \left( \int_0^1 \frac{x^{2n} dx}{1+x} \right) \left( \int_0^1 \frac{y^{2n} dy}{1+y} \right) = \sum_{n=1}^{\infty} (-1)^n \int_0^1 \int_0^1 \frac{x^{2n} y^{2n} dx dy}{(1+x)(1+y)}$$

$$= \int_0^1 \int_0^1 \left( \sum_{n=1}^{\infty} (-x^2 y^2)^n \right) \frac{dx\, dy}{(1+x)(1+y)} = -\int_0^1 \int_0^1 \frac{x^2 y^2 dx\, dy}{(1+x^2 y^2)(1+x)(1+y)} . \qquad (5)$$

Here we set $y = \frac{u}{x}$ to get

$$-\sigma = \int_0^1 \left( \int_0^x \frac{u^2 du}{(1+u^2)(u+x)} \right) \frac{dx}{(1+x)}$$

$$= \int_0^1 \left( \frac{x^2}{1+x^2} \ln 2 + \frac{\ln(1+x^2)}{2(1+x^2)} - \frac{x \arctan x}{1+x^2} \right) \frac{dx}{(1+x)}$$

$$= \ln 2 \int_0^1 \frac{x^2 dx}{(1+x^2)(1+x)} + \frac{1}{2} \int_0^1 \frac{\ln(1+x^2)}{(1+x^2)(1+x)} dx + \int_0^1 \frac{-x \arctan x}{(1+x^2)(1+x)} dx \qquad (6)$$

i.e.

$$-\sigma = A \ln 2 + \frac{1}{2} B + C \qquad (7)$$

where A, B, C are the corresponding integrals in (6). We evaluate them one by one. The first one is very easy

$$A = \frac{3}{4}\ln 2 - \frac{\pi}{8} \ . \tag{8}$$

Next,

$$B = \frac{1}{2}\left(\int_0^1 \frac{\ln(1+x^2)dx}{1+x} + \int_0^1 \frac{\ln(1+x^2)dx}{1+x^2} - \int_0^1 \frac{x\ln(1+x^2)dx}{1+x^2}\right) . \tag{9}$$

We have

$$\int_0^1 \frac{x\ln(1+x^2)\,dx}{1+x^2} = \frac{1}{2}\int_0^1 \ln(1+x^2)\,d\ln(1+x^2) = \frac{1}{4}(\ln 2)^2 \ ; \tag{10}$$

$$\int_0^1 \frac{\ln(1+x^2)\,dx}{1+x^2} = \frac{\pi}{2}\ln 2 - G , \tag{11}$$

where $G$ is Catalan's constant. See Appendix 1. Also,

$$\int_0^1 \frac{\ln(1+x^2)}{1+x}\,dx = \frac{3}{4}(\ln 2)^2 - \frac{\pi^2}{48} \ . \tag{12}$$

See Appendix 2.

Therefore,

$$B = \frac{1}{2}\left(\frac{1}{2}(\ln 2)^2 - \frac{\pi^2}{48} + \frac{\pi}{2}\ln 2 - G\right). \tag{13}$$

Finally,

$$\int_0^1 \frac{-x\arctan x}{(1+x^2)(1+x)}dx = \frac{1}{2}\int_0^1 \frac{\arctan x}{1+x}dx - \frac{1}{2}\int_0^1 \frac{x\arctan x}{1+x^2}dx - \frac{1}{2}\int_0^1 \frac{\arctan x}{1+x^2}dx , \tag{14}$$

where

$$\int_0^1 \frac{\arctan x}{1+x}dx = \frac{\pi}{8}\ln 2 , \tag{15}$$

(see Appendix 3).

$$\int_0^1 \frac{\arctan x}{1+x^2}dx = \frac{1}{2}(\arctan x)^2 \Big|_0^1 = \frac{\pi^2}{32} \ ; \tag{16}$$

$$\int_0^1 \frac{x\arctan x}{1+x^2}dx = \frac{\pi}{8}\ln 2 - \frac{1}{2}\int_0^1 \frac{\ln(1+x^2)}{1+x^2}dx = \frac{\pi}{8}\ln 2 - \frac{1}{2}\left(\frac{\pi}{2}\ln 2 - G\right) = \frac{1}{2}G - \frac{\pi}{8}\ln 2 \tag{17}$$

(after integration by parts and using (11)).

Thus

$$C = \frac{1}{2}\left(\frac{\pi}{4}\ln 2 - \frac{\pi^2}{32} - \frac{G}{2}\right). \qquad (18)$$

Now from (7) we obtain (1).

**Appendix 1**

$I_1 = \int_0^1 \frac{\ln(1+x^2)}{1+x^2} dx$ can be found in tables, see for instance 4.295.5 in [4]. The antiderivative

$\int \frac{\ln(1+x^2)}{1+x^2} dx$ was evaluated in [5, p.144] in terms of dilogarithms. We shall give now an

independent evaluation of $I_1$ based on the popular integral representation of the Catalan constant

$$G = \int_0^1 \frac{\arctan x}{x} dx = -\int_0^1 \frac{\ln x}{1+x^2} dx \ .$$

With the substitution $x = \tan\theta$,

$$I_1 + G = \int_0^1 \frac{\ln(1+x^2) - \ln x}{1+x^2} dx = -\int_0^{\pi/4} \ln\cos\theta + \ln\sin\theta \ d\theta = -\int_0^{\pi/4} \ln\frac{\sin 2\theta}{2} d\theta$$

$$= \frac{\pi \ln 2}{4} - \int_0^{\pi/4} \ln\sin 2\theta \ d\theta = \frac{\pi \ln 2}{4} - \frac{1}{2}\int_0^{\pi/2} \ln\sin\theta \ d\theta.$$

The last logsine integral is well-known (see, for example, p. 246 in [1]),

$$\int_0^{\pi/2} \ln\sin\theta \ d\theta = -\frac{\pi \ln 2}{2} \ .$$

Thus $I_1 + G = \frac{\pi \ln 2}{2}$ and $I_1 = \frac{\pi \ln 2}{2} - G$ . QED

For completeness we solve here also the logsine integral. With the substitution $\theta \to \frac{\pi}{2} - \theta$,

$$S = \int_0^{\pi/2} \ln\sin\theta \ d\theta = -\int_{\pi/2}^0 \ln\cos\theta \ d\theta = \int_0^{\pi/2} \ln\cos\theta \ d\theta,$$

Then

$$2S = \int_0^{\pi/2} \ln\sin\theta \ d\theta + \int_0^{\pi/2} \ln\cos\theta \ d\theta = \int_0^{\pi/2} \ln\frac{\sin 2\theta}{2} d\theta = \int_0^{\pi/2} \ln\sin 2\theta \ d\theta - \frac{\pi \ln 2}{2}$$

$$= \frac{1}{2}\int_0^{\pi} \ln \sin \theta \, d\theta - \frac{\pi \ln 2}{2}.$$

However,

$$\int_0^{\pi} \ln \sin \theta \, d\theta = \int_0^{\pi/2} \ln \sin \theta \, d\theta + \int_{\pi/2}^{\pi} \ln \sin \theta \, d\theta = 2\int_0^{\pi/2} \ln \sin \theta \, d\theta = 2S,$$

since the substitution $\theta \to \pi - \theta$ provides

$$\int_{\pi/2}^{\pi} \ln \sin \theta \, d\theta = \int_0^{\pi/2} \ln \sin(\pi - \theta) \, d\theta = \int_0^{\pi/2} \ln \sin \theta \, d\theta.$$

This way we come to the equation $2S = S - \dfrac{\pi \ln 2}{2}$ from which $S = -\dfrac{\pi \ln 2}{2}$.

**Appendix 2**

Proof of the evaluation $I_2 = \displaystyle\int_0^1 \frac{\ln(1+x^2)}{1+x} dx = \frac{3}{4}(\ln 2)^2 - \frac{\pi^2}{48}$.

Consider the function $F(\alpha) = \displaystyle\int_0^1 \frac{\ln(1+\alpha^2 x^2)}{1+x} dx$. Then

$$F'(\alpha) = \int_0^1 \frac{2\alpha x^2}{(1+\alpha^2 x^2)(1+x)} dx = \frac{2\alpha}{1+\alpha^2}\left[\int_0^1 \frac{dx}{1+x} + \int_0^1 \frac{x}{1+\alpha^2 x^2} dx - \int_0^1 \frac{1}{1+\alpha^2 x^2} dx\right]$$

$$= \frac{2\alpha}{1+\alpha^2}\left[\ln 2 + \frac{\ln(1+\alpha^2)}{2\alpha^2} - \frac{\arctan \alpha}{\alpha}\right] = \frac{2\alpha}{1+\alpha^2}\ln 2 + \frac{\ln(1+\alpha^2)}{\alpha(1+\alpha^2)} - \frac{2\arctan \alpha}{1+\alpha^2}.$$

What we need is

$$I_2 = F(1) = \ln 2 \int_0^1 \frac{2\alpha}{1+\alpha^2} d\alpha + \int_0^1 \frac{\ln(1+\alpha^2)}{\alpha(1+\alpha^2)} d\alpha - \int_0^1 \frac{2\arctan \alpha}{1+\alpha^2} d\alpha.$$

The first and the third integrals are trivial. With the substitution $\alpha^2 = t$ we can simplify the integral in the middle

$$\int_0^1 \frac{\ln(1+\alpha^2)}{\alpha(1+\alpha^2)} d\alpha = \frac{1}{2}\int_0^1 \frac{\ln(1+\alpha^2)}{\alpha^2(1+\alpha^2)} 2\alpha\, d\alpha = \frac{1}{2}\int_0^1 \frac{\ln(1+t)}{t(1+t)} dt,$$

and therefore,

$$I_2 = (\ln 2)^2 + \frac{1}{2}\int_0^1 \frac{\ln(1+t)}{t} dt - \frac{1}{2}\int_0^1 \frac{\ln(1+t)}{1+t} dt - \frac{\pi^2}{16}$$

$$= (\ln 2)^2 + \frac{\pi^2}{24} - \frac{1}{4}(\ln 2)^2 - \frac{\pi^2}{16} = \frac{3}{4}(\ln 2)^2 - \frac{\pi^2}{48},$$

since

$$\int_0^1 \frac{\ln(1+t)}{t} dt = \int_0^1 \sum_{n=1}^{\infty} \frac{(-1)^{n-1} t^{n-1}}{n} dt = \sum_{n=1}^{\infty} \frac{(-1)^{n-1}}{n} \int_0^1 t^{n-1} dt = \sum_{n=1}^{\infty} \frac{(-1)^{n-1}}{n^2} = \frac{\pi^2}{12}.$$

Done!

**Appendix 3**

Here we prove

$$I_3 = \int_0^1 \frac{\arctan x}{1+x} dx = \frac{\pi}{8}\ln 2.$$

This integral was evaluated in Problem 883 in [2]; it also follows from entry 4.535.1 in [1] (with $p=1$).

We shall use again differentiation on a parameter. Set $H(\alpha) = \int_0^1 \frac{\arctan(\alpha x)}{1+x} dx$. Then

$$H'(\alpha) = \int_0^1 \frac{x}{(1+\alpha^2 x^2)(1+x)} dx = \frac{-1}{1+\alpha^2}\int_0^1 \frac{dx}{1+x} + \frac{\alpha^2}{1+\alpha^2}\int_0^1 \frac{x\,dx}{1+\alpha^2 x^2} + \frac{1}{1+\alpha^2}\int_0^1 \frac{dx}{1+\alpha^2 x^2}$$

$$= \frac{-\ln 2}{1+\alpha^2} + \frac{1}{2}\frac{\ln(1+\alpha^2)}{1+\alpha^2} + \frac{\arctan \alpha}{\alpha(1+\alpha^2)} = \frac{-\ln 2}{1+\alpha^2} + \frac{1}{2}\frac{\ln(1+\alpha^2)}{1+\alpha^2} + \frac{\arctan \alpha}{\alpha} - \frac{\alpha \arctan \alpha}{1+\alpha^2}.$$

Integrating for $\alpha$ yields

$$H(1) = -\frac{\pi}{4}\ln 2 + \frac{1}{2}\int_0^1 \frac{\ln(1+\alpha^2)}{1+\alpha^2}d\alpha + \int_0^1 \frac{\arctan\alpha}{\alpha}d\alpha - \int_0^1 \frac{\alpha\arctan\alpha}{1+\alpha^2}d\alpha.$$

For the last integral we have

$$\int_0^1 \frac{\alpha\arctan\alpha}{1+\alpha^2}d\alpha = \frac{1}{2}\int_0^1 \arctan\alpha\, d\ln(1+\alpha^2) = \frac{\pi\ln 2}{8} - \frac{1}{2}\int_0^1 \frac{\ln(1+\alpha^2)}{1+\alpha^2}d\alpha,$$

so that

$$H(1) = -\frac{3\pi}{8}\ln 2 + \int_0^1 \frac{\ln(1+\alpha^2)}{1+\alpha^2}d\alpha + \int_0^1 \frac{\arctan\alpha}{\alpha}d\alpha.$$

Now according to Appendix 1, $H(1) = -\frac{3\pi}{8}\ln 2 + I_1 + G$,

and finally $I_3 = -\frac{3\pi}{8}\ln 2 + \frac{\pi}{2}\ln 2 = \frac{\pi}{8}\ln 2$.